\documentstyle{amsppt}
\magnification=\magstep1
\pageheight{23truecm}
\NoBlackBoxes

\def\Z{{\Bbb Z}}
\def\R{{\Bbb R}}
\def\Sb{{\Bbb S}}

\def\noi{\noindent}
\def\lk{\operatorname{lk}}

\def\int{\operatorname{int}}
\def\rank{\operatorname{rank}}
\def\ss{\smallskip}
\def\ms{\medskip}
\def\bs{\bigskip}

\def\kw{\mathbin{\raise 1pt \hbox {$\scriptstyle\rlap
{$\scriptstyle \sqcap $} \sqcup $}}}

 \catcode`@=11
\def\dwiestrz#1#2{\,\vcenter{\m@th\ialign{##\crcr
      $\hfil\scriptstyle{\ #1\ }\hfil$\crcr
      \noalign{\kern2\p@\nointerlineskip}
      \rightarrowfill\crcr
      \noalign{\kern\p@\nointerlineskip}
      \leftarrowfill\crcr\noalign{\kern2\p@\nointerlineskip}
      $\hfil\scriptstyle{\ #2\ }\hfil$\crcr}}\,}
\catcode`@=13

\rightheadtext{Extension of the Borsuk theorem}
\leftheadtext{J. Krasinkiewicz and S. Spie{\D z}}

\topmatter


\null\vskip1cm
\title
Extension of the Borsuk theorem on non-embeddability of spheres \vskip-3mm
\endtitle

\author
J{\' o}zef Krasinkiewicz and Stanis{\l}aw Spie{\D
z} \vskip2mm
\endauthor

\address
The Institute of Mathematics, Polish Academy of Sciences, ul.{\'
S}niadeckich 8, 00-950, Warsaw, Poland \endaddress
\email jokra$\@$impan.gov.pl
\endemail

\address
The Institute of Mathematics, Polish Academy of Sciences, ul.{\'
S}niadeckich 8, 00-950, Warsaw, Poland
\endaddress
\email spiez$\@$impan.gov.pl
\endemail

\abstract It is proved that the suspension
$\sum M$ of a closed $n$-dimensional manifold $M$, $n\ge1$,
does not embed in a product of $n+1$ curves. In fact,
the ultimate result will be proved in a much more general setting.
This is a far-reaching generalization the Borsuk theorem on
non-embeddability of the sphere $\Sb^{n+1}$ in a product of $n+1$ curves.
\vskip4mm

\bs\centerline{\bf CONTENTS}

\ms \item{\eightrm 1.} {\eightit Introduction \hfill{2}}\ss

\item{\eightrm 2.} {\eightit On
quasi-manifolds, weak manifolds and
ramified manifolds \hfill{2}}\ss

\item{\eightrm 3.} {\eightit Second Factorization Theorem
\hfill{4}}\ss

\item{\eightrm 4.} {\eightit  Proof of Theorem 1.2 \hfill{5}}\ss

\item{\eightrm 5.} {\eightit  Factorization Theorem for ramified 3-manifolds \hfill{5}}\ss

\item{\eightrm 6.} {\eightit  Product structure of locally connected
weak manifolds lying in products \hfill \vskip-1mm of curves \hfill{6}}\ss

\item{\eightrm 7.} {\eightit  Proof of Theorem 1.3 \hfill{6}}\ss
\item{} {\eightit Problem \hfill{7}}

\item{}{\eightit References \hfill{7}}
\endabstract



\endtopmatter

\vfill \break

\document

\centerline{1. {\it Introduction}}\bs

All {\it spaces} discussed in this note are metrizable and all
{\it mappings} (also called {\it maps}) are continuous. By a {\it
compactum} we mean a compact metric space, by a {\it continuum} --
connected compactum, and by a {\it curve} -- 1-dimensional
continuum.

In 1975 K.~Borsuk \cite{B} discovered a remarkable property of
spheres (this property provided an answer to a problem of
J. Nagata \cite{N, p. 163}):

\proclaim{Theorem 1.1 (Borsuk)} The $n$-sphere $\Sb^n$, $n\ge 2$,
does not embed in a product of n curves.\qed
\endproclaim

In this article we present a simple elementary proof of
the following extension of this theorem.

\proclaim{Theorem 1.2} The suspension $\sum M$ of any
closed $n$-manifold $M$, where $n\ge1$, does not embed in a product of $n+1$ curves.
\qed\endproclaim

The proof is given in Section 4. Moreover, we shall prove the following more
general result (see Section 6 for a proof).

\proclaim{Theorem 1.3} The suspension $\sum X$ of any locally
connected quasi $n$-manifold $X$, where $n\ge1$, does
not embed in a product of $n+1$ curves.
\qed\endproclaim

The proof is based on the Second Factorization Theorem 3.2
and the Second Structure Theorem 6.1 that are stronger forms of
the corresponding results from \cite{K-K-S}.
In addition we show that the Second Factorization Theorem holds for
ramified 3-manifolds, see Theorem 5.3.\bs

\centerline{2. {\it On quasi-manifolds, weak manifolds and ramified manifolds}}\bs

An $n$-dimensional compactum is
said to be an {\it $n$-manifold at a point} $x$ if there
is an open neighborhood $W$ of $x$ in $X$ that is an open $n$-disc.
Obviously, $X$ is an $n$-manifold at every point $x\in X$ if and
only if $X$ is a closed $n$-manifold. The property discussed in the next
paragraph follows from the Borsuk Separation
Theorem which relates compact sets separating ${\Sb}^n$, $n\ge1$, to
essential mappings from those sets into ${\Sb}^{n-1}$
(cf. \cite{E-S, p. 302}).

For an $n$-dimensional compactum $X$ which is an $n$-manifold at a point
$x_0 \in X$ there is an open neighborhood $W$ of $x_0$ in $X$ such
that for every open set $U\ni x$ with $\overline U\subset W$ there is
an essential map $\partial U\to \Sb^{n-1}$.
$($In fact, this holds for every neighborhood $W$ that is an open
$n$-disc.$)$ One can use this property to define a larger
class of $n$-dimensional compacta.

Namely, an $n$-dimensional compactum $X$ is said to be
a {\it quasi} $n$-{\it manifold at a point} $x \in X$ if there is
an open neighborhood $W$ of $x$ in $X$ such that for every open set
$U\ni x$ with $\overline U\subset W$ and $\dim \partial U\le n-1$ there is
an essential map $\partial U\to \Sb^{n-1}$; equivalently,
$\widetilde H^{n-1}(\partial U)\ne0$. (By $H^{\ast}(\cdot)$ we denote
the {\v Cech} cohomology with integer coefficients $\Z$.)
If $X$ is a quasi
$n$-manifold at $x$ then every neighborhood of $x$ is
$n$-dimensional. A compactum $X$ is said to be a {\it quasi
$n$-manifold} if it is a quasi $n$-manifold at every point. For polyhedra
we have the following characterization.

\proclaim{Theorem 2.1} (see {\cite{K-K-S, Theorem 2.2}{\rm )} 
For a compact $n$-dimensional
polyhedron $X=|K|$ and a point $x\in X$ the following are equivalent:\ss

{\rm (i)} $X$ is a quasi $n$-manifold at $x$;\ss

{\rm (ii)} $\widetilde H^{n-1}(|\lk(x,L)|)\ne0$
for every triangulation $L$ of $X$ with $x$ as a vertex;\ss

{\rm (iii)} either $\dim\sigma(x)=n$ or $\dim\sigma(x)<n$ and
$\widetilde H^{n-\dim\sigma(x)-1}(|\lk(\sigma(x),K)|)\ne0$,
where $\sigma(x)$ is the carrier of $x$ in $K$;\ss

{\rm (iv)} $H^n(X,X\setminus \{x\})\neq 0$. \qed\endproclaim

Let $X$ be a compact polyhedron. A triangulation $K$ of $X$ is
said to be {\it fine} if for each point $x\in X$ there is a vertex
$v$ of $K$ such that $(X,x)$ is homeomorphic to $(X,v)$, written
$(X,x)\approx (X,v)$. Every polyhedron has a fine triangulation,
e.g., the barycentric subdivision of a triangulation. We have also
the following description of quasi $n$-manifolds among $n$-dimensional
polyhedra. The proof is similar to the proof of the above theorem.

\proclaim{Theorem 2.2} For a compact $n$-dimensional polyhedron
$X$ the following are equivalent:\ss

{\rm (i)} $X$ is a quasi $n$-manifold;\ss

{\rm (ii)} there is a fine triangulation of $X$ such that the
links of vertices have nontrivial $(n-1)$-cohomology
$H^{n-1}(.)$.\qed
\endproclaim

A compact pair $(X,A)$ is said to {\it extend} to a pair $(X',A')$ if
there is an embedding $h:X\to X'$ such that $h(A)\subset A'$.

A compact $n$-dimensional space $X$ is said to be a
{\it weak n-manifold at a point} $x\in X$ if there is
an open neighborhood $W$ of $x$ in $X$ such that no pair
$(\overline{U},\partial U)$, where $U\ni x$ is open and
$\overline{U}\subset W$, extends to a compact pair $(P,Q)$ of absolute
retracts with $\dim (P,Q)=(n,n-1)$ (that is, $\dim P=n$ and
$\dim Q=n-1$). An $n$-dimensional compactum is said to be a
{\it weak $n$-manifold} if it is a weak $n$-manifold at every point.
If $X$ is a quasi $n$-manifold at a point $x$ then $X$ is a weak
$n$-manifold at $x$ as well.

Let $S$ be a subset of a compactum $X$. Then $X$ is said to be an
{\it $n$-manifold} ({\it quasi $n$-manifold, weak $n$-manifold,
{\rm resp.}}) {\it off $S$} if $X$ is $n$-dimensional at every
point and $X$ is an $n$-manifold (quasi $n$-manifold, weak $n$-manifold,
resp.) at each point $x\in X\setminus S$.

A cell complex $\,$ \footnote{The notion of {\it cell complex}
and other related notions are understood in the sense of \cite{R-S}.
By $K^{(m)}$ we denote the $m$th skeleton
of $K$; if $m<0$ then $K^{(m)}=\emptyset$ by definition.}
$K$ is called a {\it ramified
$n$-complex} if $K$ consists of a finite collection of $n$-cells
and their faces such that each $(n-1)$-cell is a face of at least
two $n$-cells. (A ramified 0-complex consist of at least two points.)
The underlying polyhedron $X=|K|$ of a ramified $n$-complex $K$
is called a {\it ramified $n$-manifold}. Notice that a ramified
$n$-manifold $X=|K|$ is a quasi
$n$-manifold off $|K^{(n-2)}|$. Any polyhedron
that is a quasi $n$-manifold is a ramified $n$-manifold.
For $n\le2$ the converse
is also true. \footnote{For $n=3$ it fails (see \cite{K-K-S}).}

\proclaim{Lemma 2.3} Let $X$ be a weak $n$-manifold off a set with
dimension $\le n-2$. If $U$ is a non-void open subset of X and
$h:U\to \R^n$ is an embedding such that $h(U)$ is closed in $\R^n$ then $h(U)=\R^n$.
\endproclaim

\demo{Proof} Suppose $h(U)$ is a proper subset of $\R^n$. Assume $X$ is a weak
$n$-manifold off a set $S$ with $\dim S\le n-2$.
There is a point $a\in h(U)\setminus h(S)$ because $X$ is $n$-dimensional
at each point. Pick any point $b\in \R^n \setminus h(U)$. Since
$\dim h(S) \le n-2$, the Mazurkiewicz theorem \cite{E, Theorem 1.8.18, p. 62}
implies that there
is a continuum $C\subset \R^n \setminus h(S)$ connecting $a$ with $b$. Let
$C_0$ be the component of $C\setminus h(U)$ containing $b$. It follows that
$\overline{C_0}\cap \partial h(U)\neq \emptyset$ \cite{Kur, p. 172}.
Pick any point $c\in \overline{C_0}\cap\partial h(U)$. Then
$h^{-1}(c)\in U\setminus S$. Let $B_r$ denote the
open $n$-ball in $\R^n$ with radius $r>0$ and center at $c$, and let $S_r$
denote its boundary. Since $X$ is a weak $n$-manifold at $h^{-1}(c)$ there
is an open
neighborhood $W$ of $h^{-1}(c)$ in $X$ satisfying the definition of weak
$n$-manifold at $h^{-1}(c)$.
Then there is a positive number $r<d(c,b)$ such that
$\overline{V}\subset h(U\cap W)$
for $V=B_r\cap h(U)$. Then $\partial_{h(U)}V\subset S_r\setminus C_0$.
Since $C_0\cap S_r \neq\emptyset$ the pair  $(\overline{V},\partial_{h(U)}V)$
extends to a pair of absolute retracts
$(\overline{B_r},D)$ with dimension $(n,n-1)$, where $D$ is a closed
$(n-1)$-disc in $S_r$. Since $h^{-1}(V)$ is an open neighborhood of
$h^{-1}(c)$ whose closure lies in $W$, and $h$ transforms
$(\overline{h^{-1}(V)},\partial h^{-1}(V))$ homeomorphically onto
$(\overline{V},\partial_{h(U)}V)$ we get a contradiction.
\enddemo

\proclaim{Theorem 2.4} Let $X$ be a weak $n$-manifold off a set with
dimension $\le n-2$. If $f:X\to |K|$
is an embedding, where $K$ is an $n$-dimensional cell complex,
then $f(X)=|L|$, where $L$ is a ramified $n$-subcomplex of $K$.
In particular, the conclusion holds for ramified $n$-manifolds.
\endproclaim

The proof is based on the preceding lemma and runs the same way as the proof
of Theorem 2.7 in \cite{K-K-S}.\qed

\bs\centerline{3. {\it Second Factorization Theorem}}
\bs
In \cite{K-K-S} the following theorem has been established.

\proclaim{Factorization Theorem 3.1} Let $X$ be a locally
connected weak $n$-manifold with $H^1(X)$ of finite rank. If $f =
(f_1, \cdots ,f_n):X \to Y_1 \times \cdots \times Y_n$ is an
embedding in a product of $n$ curves, then there exist mappings $g
= (g_1, \cdots ,g_n):X \to P_1 \times \cdots \times P_n$ and $h=
h_1 \times \cdots \times h_n:P_1 \times \cdots \times P_n \to Y_1
\times \cdots \times Y_n$ such that $f_i = h_i \circ g_i$ for each
$i=1,\cdots ,n$ $($hence $f=h \circ g$$)$, where $g_i:X \to P_i$
is a monotone surjection, $P_i$ is a graph with no endpoint, and
$h_i:P_i\to Y_i$ is $0$-dimensional. In particular, $g$ is an embedding.
\qed\endproclaim

Let us recall the main steps in the proof of this theorem. The existence of
the mappings follows from the Whyburn factorization theorem \cite{W, p.
141} (which is applicable to any compactum). The local connectedness of $X$
plus $\rank H^1(X)<\infty$
ensure that $P_i$ is a compact 1-dimensional ANR with\ss
    (*) $H^1(P_i)$ {\it isomorphic to a subgroup of} $H^1(X)$.\ss
\noi Then the following basic observation was made:\ss
    (**) {\it if $g_i(x)$ is an endpoint of $P_i$ then $X$ is not a weak
    $n$-manifold at $x$}.\ss
\noi Let $NWM(X)$ designate the subset of $X$ composed of all points at which $X$
is not a weak $n$-manifold, and let $E(P_i)$ designate the set of endpoints
of $P_i$. From (**) it follows that\ss
    (***) $\# E(P_i) \le \# NWM(X)$.\ss
\noi One easily sees that if $P$ is a compact 1-dimensional ANR and the set $E(P)$
of its endpoints is finite then $P$ is a graph. Therefore, by a slight
modification of the discussed proof, one obtains the following refined version
of the Factorization Theorem:

\proclaim{Second Factorization Theorem 3.2} Let $X$ be a locally
connected weak $n$-manifold off a finite set with $H^1(X)$ of finite
rank. If $f =
(f_1, \cdots ,f_n):X \to Y_1 \times \cdots \times Y_n$ is an
embedding in a product of $n$ curves, then there exist mappings $g
= (g_1, \cdots ,g_n):X \to P_1 \times \cdots \times P_n$ and $h=
h_1 \times \cdots \times h_n:P_1 \times \cdots \times P_n \to Y_1
\times \cdots \times Y_n$ such that $f_i = h_i \circ g_i$ for each
$i=1,\cdots ,n$, where $g_i:X \to P_i$ is a monotone surjection onto
a graph satisfying {\rm (**)}, and
$h_i:P_i\to Y_i$ is $0$-dimensional. Consequently,
$\# E(P_i)\le \# NWM(X)$ for each $i$, and $g$ is an embedding.
\qed\endproclaim\bs

\centerline{4. {\it Proof of Theorem 1.2}}\bs

Suppose $\sum M$ embeds in a product of $(n+1)$ curves. We may assume
that $M$ is connected.
Then, by the Second Factorization Theorem 3.2, there is a mapping $g
= (g_1, \cdots ,g_{n+1}):\sum M\to P_1 \times \cdots \times P_{n+1}$
such that each $g_i:\sum M \to P_i$ is a monotone surjection with
the following properties: $P_i$ is a compact connected
1-dimensional ANR, $H^1(P_i)=0$, and $P_i$ has at most two endpoints
(because $\sum M$ is an $(n+1)$-manifold off the vertices).
It follows that $P_i$ is an arc
(because it is a dendrite with at most two endpoints). Let $K$ denote the
natural cell structure on $P_1 \times \cdots \times P_{n+1}$ with just
one $(n+1)$-cell. Since $g:\sum M\to |K|$ is an embedding and $\sum M$ is
an $(n+1)$-manifold off the vertices,
$g(\sum M)=P_1 \times \cdots \times P_{n+1}$ by Theorem 2.4. Thus the
$(n+1)$-cell $P_1 \times \cdots \times P_{n+1}$ is an $(n+1)$-manifold
off two points, a contradiction.\qed\ms

\noi{\bf Remark.} Notice that Theorem 2.1 also follows from \cite{D-K}, but
the above proof is more elementary.\bs

\centerline{5. {\it Factorization Theorem for ramified 3-manifolds
}}\bs

\proclaim{Lemma 5.1} If $X=|K|$ is a ramified $n$-manifold with a
triangulation $K$ then for every simplex $\sigma\in K$ with
$\dim \sigma<n$ the link $\lk(\sigma,K)$ is a ramified
$(n-\dim\sigma -1)$-complex. Also, $X$ is a
quasi $n$-manifold off $|K^{(n-3)}|$ $($hence a fortiori a weak
$n$-manifold off this set$)$.
\endproclaim

\demo{Proof} The proof of the first assertion is straightforward.
To prove the second one, first notice that $X$ is a quasi $n$-manifold
off $|K^{(n-2)}|$ (because $X\setminus |K^{(n-2)}|$ is a union
of open $n$-cells). Next consider the case where
$x\in |K^{(n-2)}|\setminus |K^{(n-3)}|$.
Then $\dim \sigma(x)=n-2$, where $\sigma(x)$ is the carrier of $x$ in $K$.
Hence $|\lk(\sigma(x),K)|$ is a ramified 1-manifold. Consequently,
$H^1(|\lk(\sigma(x),K)|)\neq0$. The conclusion follows now from
Theorem 2.1(iii).\qed
\enddemo

\proclaim{Corollary 5.2} Let $X=|K|$ be a ramified $3$-manifold with a
triangulation $K$. Then $X$ is a quasi $3$-manifold off the vertices of $K$.
\qed\endproclaim

Combining this corollary with the Second Factorization Theorem 3.2 we
obtain the following

\proclaim{Theorem 5.3} If a ramified $3$-manifold embeds in a
product of three curves then it embeds in
a product of three graphs as well. Moreover, if
$f =(f_1,f_2,f_3):X \to Y_1 \times Y_2\times Y_3$
is an embedding in a product of curves then there exist mappings $g
= (g_1,g_2,g_3):X \to P_1 \times P_2\times P_3$ and $h=
h_1 \times h_2 \times h_3:P_1 \times P_2\times P_3 \to Y_1
\times Y_2 \times Y_3$ such that $f_i = h_i \circ g_i$ for each
$i=1,2,3$, where $g_i:X \to P_i$
is a monotone surjection onto a graph, and
$h_i:P_i\to Y_i$ is $0$-dimensional. Then $g$ is an embedding.
\qed\endproclaim\bs

\centerline{6. {\it Product structure of locally connected weak
manifolds}} \centerline{{\it lying in products of curves}}\bs

Using the Second Factorization Theorem in the proof of the Structure
Theorem 5.1 in \cite{K-K-S} we obtain the following
improvement of the Structure Theorem. Denote
$$J_X=\{j\in\{1,\cdots,n\}:pr_j(X)~ \text{\it {is a circle}}\}.$$

\proclaim{Second Structure Theorem 6.1} Let $X$ be a connected and
locally connected weak $n$-manifold off a finite set, 
with $H^1(X)$ of finite rank, lying
in a product $Y_1\times\cdots\times{}Y_n$ of $n$ curves, $n\ge1$.
Then $X$ is a polyhedron and the following conditions are
fulfilled: \ss{\rm(1)} $\rank{}H^1(X)\ge n$ $;$ \ss{\rm(2)} If
$\rank{}H^1(X)=n+k$, where $k<n$, then $J_X$ contains at least
$n-k$ elements. In particular, if $\rank{}H^1(X)=n$ then
$J_X=\{1,\cdots,n\}$ $;$ \ss{\rm(3)} If $J_X=\{1,\cdots,n\}$ then
$X=pr_1(X)\times\cdots\times pr_n(X)$ is an $n$-torus $;$
\ss{\rm(4)} If $J_X$ is a proper non-void subset of
$\{1,\cdots,n\}$ then $X=(\prod_{j\in J_X}pr_j(X))\times
pr_{J_X^c}(X)$, where the first factor is an $n_{J_X}$-torus, and
$pr_{J_X^c}(X)$ is a polyhedron that is a weak
$n_{J_X^c}$-manifold in $\prod_{j\in J_X^c}Y_j$ having no
projection onto a circle.

If, in addition, $X$ is a quasi $n$-manifold, then $pr_{J_X^c}(X)$
is a quasi $n_{J_M^c}$-manifold.
\qed\endproclaim\bs

\centerline{7. {\it Proof of Theorem 1.3}}\bs

Let $X$ be a locally connected quasi $n$-manifold, $n\ge1$,
and suppose the suspension $\sum X$ embeds in a product of $n+1$ curves.
Without loss of generality we can assume that $X$ is connected.
By Corollary 2.3 in \cite{K-K-S} the suspension is a
quasi $(n+1)$-manifold off the vertices.
Hence $\sum X$ is a locally connected, weak $(n+1)$-manifold off
the vertices, and $H^1(\Sigma X)=0$. Then, by the
Second Structure Theorem 6.1(1), we reach a contradiction.\qed\bs

We end this paper with the following

\proclaim{Problem} Is the Factorization Theorem true
for ramified $n$-manifolds?
\endproclaim


\Refs \widestnumber\key{K-K-S}

\ref\key B \by K. Borsuk \paper Remarks on the Cartesian product of two
$1$-dimensional spaces \jour Bull. Polish Acad. Sci. Math. \yr
1975  \vol 23 \pages 971--973
\endref

\ref\key D-K \by J.~Dydak and A.~Koyama \paper Compacta not
embeddable into Cartesian products of curves \jour Bull. Pol.
Acad. Sci. Ser. Math. \yr 2000 \vol 48 \pages 51--56
\endref

\ref\key E-S \by S. Eilenberg and N. Steenrod \book Foundations of
Algebraic Topology \publ Princeton University Press, Princeton,
New York\yr 1952
\endref

\ref\key E \by R. Engelking \book Theory of dimensions finite and infinite \publ
Sigma Series in Pure Mathematics, 10; Heldermann Verlag, Lemgo\yr 1995
\endref

\ref\key H-W \by W. Hurewicz and H. Wallman \book Dimension Theory
\publ Princeton \yr 1948
\endref

\ref\key K-K-S \by A. Koyama, J. Krasinkiewicz, S. Spie\D z \paper
Embedding compacta into products of curves \jour Trans. Amer. Math. Soc. (submitted)
\endref


\ref\key Kur \by K.~Kuratowski \book Topology, vol. II \publ
PWN-Academic Press, Warsaw-New York\yr 1968
\endref

\ref\key N \by J. Nagata \book Modern Dimension Theory \publ
North-Holland, Amsterdam \yr 1965
\endref

\ref\key R-S \by C.~P.~Rourke and B.~J.~Sanderson \book
Introduction to Piecewise-Linear Topology \yr 1982 \publ
Springer-Verlag, Berlin, Heilderberg, New York
\endref

\ref\key Sp \by E.~Spanier \book Algebraic Topology \yr 1966 \publ
McGraw-Hill, New York
\endref

\ref\key W \by G.~T.~Whyburn \book Analytic Topology \publ Amer.
Math. Soc. Colloquium Publications \vol 28 \yr 1942
\endref

\endRefs
\enddocument